%

\magnification\magstep1
\openup 2\jot
\def\proof{{\bf Proof. }}
\def\endproof{\par\medskip}

\def\set#1{\{\,#1\,\}}
\def\seq#1{\langle\,#1\,\rangle}
\def\forces{\mathrel\|\joinrel\relbar} 
\def\rest{|}

\def\sse{\subseteq}
\def\a{\alpha}
\def\b{\beta}
\def\d{\delta}

\def\g{\gamma}
\def\h{\eta}
\def\k{\kappa}
\def\l{\lambda}
\def\m{\mu}

\def\r{\rho}
\def\s{\sigma}
\def\t{\tau}
\def\w{\omega}
\def\x{\xi}
\def\z{\zeta}
\def\A{\aleph}  
\def\cf{\mathop{\rm cf}\nolimits}

\def\range{\mathop{\rm range}\nolimits}

\outer\def\case #1. {\medbreak\noindent{\sl Case #1.\quad}}


\centerline{Extensions of the Erd\H os-Rado Theorem}
\centerline{by}
\centerline{James E. Baumgartner\footnote{$^1$}{The preparation of
this paper was partially supported by National Science Foundation grant
number
DMS--8906946.}, Andr\'as Hajnal\footnote{$^2$}{Research supported by
Hungarian
National Science Foundation OTKA grant 1908.}, and Stevo Todor\v
cevi\'c\footnote{$^3$} {Research supported by the Science Fund of Serbia
grant
number 0401A}}

\bigskip
\centerline{\it Dedicated to Eric Milner on the occasion of his third coming
of
age}

\bigskip \noindent{\bf 0. Statement of results.}  We find here some
extensions
of the Erd\H os-Rado Theorem that answer some longstanding problems. 
Ordinary
partition relations for cardinal numbers are fairly well understood (see
[EHMR]), but for ordinal numbers much has been open, and much remains open. 
For
example, any proof of the simplest version of the Erd\H os-Rado Theorem seems
to yield
$$\hbox{ For any regular cardinal $\k$, if $\m<\k$ then }
(2^{<\k})^+ \to (\k+1)^2_\m,$$
but to replace $\k+1$ by $\k+2$ seems quite a nontrivial problem.  In this
paper we will prove that if $\k$ is regular and uncountable, then

\proclaim Theorem 3.1. $\forall k<\w$ $\forall\x<\log\k$ $(2^{<\k})^+ \to
(\k+\x)^2_k$.

\proclaim Theorem 4.1. $\forall n,k<\w$ $(2^{<\k})^+ \to (\r,(\k+n)_k)^2$,
where
$\r=\k^{\w+2}+1$.

\proclaim Theorem 5.1. $\forall n<\w$ $(2^{<\k})^+\cdot\w \to (\k\cdot
n)^2_2$.

The actual version of 5.1 is slightly stronger.  Here $\log\k$ is the least
cardinal $\m$ such that $2^\m\ge\k$, the exponentiation in 4.1 is ordinal
exponentiation, and the products in 5.1 represent ordinal multiplication. 
For
$\k=\w$, 3.1, 4.1, and 5.1 all follow from the known result $\w_1\to(\a)^2_k$
for all $\a<\w_1$ and all $k<\w$ (see 
), so the uncountable case is the
interesting one.

 The
proofs make use of elementary substructures of structures of the form
$(H(\m),\in)$ where $H(\m)$ is the set of all sets hereditarily of
cardinality
$<\m$, and of ideals on ordinals generated from such elementary
substructures.  It is possible to recast our arguments in such a way as to
fit
them under the heading of ramification arguments but we choose not to do so,
both
because the proofs remain clearer this way, and because we believe that this
may generally be a better approach to ramification arguments.

The proof of 4.1 also uses a metamathematical trick seen in [BH] for
countable
ordinals.  Our approach is to prove 4.1 in a generic extension of the
universe
obtained by $\k$-closed forcing, and then to argue that it must therefore be
true.

Several open problems remain.  The simplest versions are whether any of the
following are provable from GCH:
$$\w_3\to(\w_2+2)^2_\w$$
$$\w_3\to(\w_2+\w_1,\w_2+\w)^2$$
$$\w_2\to(\w_1^{\w+2}+2,\w_1+2)^2$$
or whether the following is provable from CH:
$$\w_3\to(\w_1+\w)^2_3.$$

\bigskip

\noindent{\bf1. Terminology.} Our partition calculus notation is standard. 
If
$x$ is a set and $\l$ is a cardinal then $[x]^\l=\set{y\sse x:|y|=\l}$,
$[x]^{<\l} = \bigcup\set{[x]^\k:\k<\l}$, and $[x]^{\le\l}=[x]^{<\l^+}$. 
Usually, but not always, such a $\l$ is finite.  The ordinary (ordinal)
partition relation $\a\to(\b)^2_\m$ means that $\forall f:[\a]^2\to\m$
$\exists
X\sse\a$ $X$ has order type $\b$ and is homogeneous for $f$, i.e., $f$ is
constant on $[X]^2$.  Here $\a$ and $\b$ may be ordinals; $\m$ is always
taken
to be a cardinal.  If $f(x)=i$ for all $x\in[X]^2$ then we say $X$ is
$i$-{\it homogeneous}.  The unbalanced relation
$$\l\to(\g_0,\g_1,\ldots,\g_{n-1})^2$$
means that $\forall f:[\l]^2\to n$ $\exists X\sse\l$ $\exists i<n$ $X$ has
order type $\g_i$ and is $i$-homogeneous.  We abbreviate
$$\l\to(\a,\b,\b,\b,\ldots,\b)^2$$
by $\l\to(\a,(\b)_n)^2$ if there are $n$ occurrences of $\b$.

If $\k$ and $\l$ are cardinals, then $\l^{<\k}=\sum\set{\l^\m:\m<\k}$. 
Clearly
$2^{<\k}\ge\k$.  If $\k$ is regular then $\k^{<\k}=2^{<\k}$.

If $(N,\in)$ is an elementary substructure of $(H(\l),\in)$ (i.e.,
$(N,\in)\prec (H(\l),\in)$) then we often abbreviate this by $N\prec H(\l)$. 
The $\in$-relation is always understood.  Suppose $\k$ is regular, $N\prec
H(\k^{++})$, $|N|=\k$ and $\a=N\cap\k^+<\k^+$.  Then we may define an ideal
$I$
on $\a$ as $\set{X\sse\a:\exists A\in N\ \a\notin A \hbox{\ and }X\sse N}$. 
Clearly $I$ contains all bounded subsets of $\a$.  If in addition
$2^{<\k}=\k$
then it is quite possible that $[N]^{<\k}\sse N$, and we generally work with
sets $N$ with this property.  Note in that case that $I$ is $\k$-closed,
i.e.,
if $\m<\k$ and $X_i\in I$ for $i<\m$, then $\bigcup\set{X_i:i<\m}\in I$. 
This
follows from the fact that if $A_i\in N$ for each $i<\m$ and $\a\notin A_i$,
then
$\seq{A_i:i<\m}\in N$ since $[N]^{<\k}\sse N$ so $A=\bigcup\set{A_i:i<\m}\in
N$
and $\a\notin A$.

If $I$ is an ideal on $\a$ then $I^+=\set{X\sse\a:X\notin I}$ and
$I^*=\set{X\sse\a:\a-X\in I}$.  Note that if $I$ is not proper, then $I^+=0$
and $I^*=I$.  If $I(\a,\b)$ is an ideal then we write $I^+(\a,\b)$ instead of
$I(\a,\b)^+$.

If $N\prec H(\k^{++})$ as above, $S\in N$ and $\a\in S$, then $S\cap\k^+$
must
be cofinal in $\k^+$.  If not, an upper bound is definable from $S$ and
$\k^+$,
both of which lie in $N$ ($\k^+\in N$ since it is the largest cardinal in
$N$),
so the upper bound must be $<\a$, which is impossible.  It follows that
$S\cap
\a$ is cofinal in $\a$.  In fact, $S\cap\k^+$ is actually stationary in
$\k^+$,
since otherwise there would be closed unbounded $C\in N$ with $S\cap C=0$. 
But
then $C\cap\a$ is unbounded in $\a$ so $\a\in C$, which is impossible.

We often use variations of the ideal $I$ to produce large homogeneous sets
for
a partition function $f\in N$.  Such arguments can also be made using
ramification ideas.  A very good introduction to the ideas we plan to use is
in
Section 2, where a simple version of the Erd\H os-Rado Theorem is proved with
these techniques.

In Section 4 we will use the notions of stationary and closed unbounded
subsets
of $[X]^\k$, where $X$ is a set such that $|X|>\k^+$.  We assume the reader
is
familiar with these notions.

\bigskip

\noindent{\bf2. The Erd\H os-Rado Theorem.} In this section we present a
slightly nonstandard proof of the Erd\H os-Rado Theorem, and an associated
result, using elementary substructures.  Our goal is to prepare the way for
an
extension of the Erd\H os-Rado Theorem to be proved in the next section.

\proclaim Theorem 2.1. (Erd\H os-Rado [ER]) Let $\k$ be a regular cardinal
and
let $\l=(2^{<\k})^+$.  Then $\l\to(\k+1)^2_\m$ for all $\m<\k$.

\proof Fix $f:[\l]^2\to\m<\k$.  Choose $N\prec H(\l^+)$ such that $f\in N$,
$|N|=2^{<\k}$, $[N]^{<\k}\sse N$ and $N\cap\l=\a<\l$.  (Note that this
implies $\cf\a=\k$.  If $\cf\a<\k$ then some cofinal subset $x$ of $\a$
would be an element of $N$, and $\a=\bigcup x$ so $\a\in N$, contradiction.) 
Let
$I_\a=\set{X\sse\a: \exists A\in N\ \a\notin A \hbox{\ and }X\sse A}$.  Then
$I_\a$ is a $\k$-complete ideal on $\a$ since $[N]^{<\k}\sse N$.

For $i<\m$ let $H_i=\set{\b<\a:f\{\b,\a\}=i}$.  Then
$\bigcup\set{H_i:i<\m}=\a$ so for some $i$ $H_i\in I_\a^+$.

\proclaim Claim. If $H_i\in I_\a^+$, $X\sse H_i$ and $X\in I_\a^+$ then
$\exists Y\sse X$ $Y$ is $i$-homogeneous of order type~$\k$.

This will suffice since then $Y\cup\{\a\}$ is $i$-homogeneous of type $\k+1$.
 
It will suffice to show

\proclaim Lemma 2.2. If $X\sse H_i$, $X\in I_\a^+$, $x\sse X$, $|x|<\k$ and
$x$
is $i$-homogeneous, then $\exists\b\in X$ $\b>\sup x$ and $x\cup\{\b\}$ is
$i$-homogeneous.

\proof Let $Z=\set{\g<\l:\forall\d\in x f\{\d,\g\}=i}$.  Then $Z\in N$ and
$\a\in Z$ so $\a-Z\in I_\a$.  Since $X\in I_\a^+$ we have $X\cap Z\in I_\a^+$
as well.  Now choose $\b\in X\cap Z$, $\b>\sup x$.
\medskip
Using the same methods we can do a little better.

\proclaim Theorem 2.3. As in Theorem 2.1, let $\k$ be regular,
$\l=(2^{<\k})^+$, $\m<\k$ and $f:[\l]^2\to\m$.  Then $\exists A\sse\l$
$\exists
i<\m$ $A$ is $i$-homogeneous and either $|A|=\l$ and $i=0$ or else $A$ has
type
$\k+1$ and $i>0$.  (In short notation $\l\to(\l,(\k+1)_{\m-1})^2$.)

\proof We use the notation of the previous proof.  If $H_i\in I_\a^+$ for
some
$i>0$ we are done, so assume $H_i\in I_\a$ for all $i>0$.  Then
$\bigcup\set{H_i:i>0}\in I_\a$ so for some $X\in N$ we have $\a\in X$ and
$X\cap\a\sse H_0$.  If in $N$ we define $Y=\set{\b\in X:\forall\g\in X\cap\b\
f\{\g,\b\}=0}$, then $\a\in Y$.  Hence $|Y|=\l$ and is clearly 0-homogeneous.
\medskip

\noindent {\bf Remark.} Theorem 2.3 is a well known result of Erd\H os,
Dushnik and Miller.  See [DM], for example.  To get the rest of the Erd\H
os-Rado
Theorem with this approach (i.e., the version for partitions of $n$-element
sets) simply proceed via end-homogeneous sets, an approach which is itself
highly
adaptable---and well known---in the context of elementary substructures. 
Details are left to the reader.

\bigskip

\noindent{\bf 3. A balanced extension of the Erd\H os-Rado Theorem.}  In this
section we seek an improvement of the Erd\H os-Rado Theorem, Theorem 2.1. 
This
will require a clearer analysis of the ideals $I_\a$ defined in Section 2.

Recall that $\log\k$ is the least cardinal $\m$ such that $2^\m\ge\k$.

\proclaim Theorem 3.1. Let $\k$ be a regular uncountable cardinal, and let
$\l=(2^{<\k})^+$.  Then $\forall k<\w$ $\forall\x<\log\k$\ \ 
$\l\to(\k+\x)^2_k$.

Some parts of this theorem are already known.  In unpublished work, Hajnal
proved Theorem 3.1 for $\k=\w_1$ and $k=2$ about thirty years ago.  Shelah
proved
Theorem 3.1 for the case $k=2$ in [Sh1, \S6].

The rest of this section will be devoted to the proof of Theorem 3.1.

Fix $f:[\l]^2\to k$, and let $\seq{N_\a:\a<\l}$ be a continuous sequence of
elementary substructures of $H(\l^+)$ such that $f\in N_0$, $|N_\a|=2^{<\k}$,
$N_\a\in N_{\a+1}$ and $[N_\a]^{<\k}\sse N_{\a+1}$.  We call a sequence
continuous if each limit point is the union of the preceding elements.  We
may
assume that $\set{\a:N_\a\cap\l=\a}$ is closed and unbounded in $\l$.  Note
that
if $\a\in S_0=\set{\a:N_\a\cap\l=\a \hbox{\ and }\cf\a=\k}$ then
$[N_\a]^{<\k}\sse N_\a$.

If $\a\in S_0$ then let $I_\a=\set{X\sse\a:\exists A\in N_\a\ \a\notin A
\hbox{\ and }X\sse A}$.

\proclaim Lemma 3.2. Let $S\sse S_0$ be stationary.  Then there is closed
unbounded $C\sse\l$ such that $\forall\a\in C\cap S$ $S\cap\a\in I_\a^+$.

We refer to $\a$ as a {\it reflection point\/} of $S$.

\noindent{\bf Proof of Lemma 3.2.} Let $\seq{M_\a:\a<\l}$ be a continuous
sequence of elementary substructures of $H(\l^+)$ such that $N_\a\sse M_\a$,
$|M_\a|=2^{<\k}$ and $S\in M_0$.  Let $N=\bigcup\set{N_\a:\a<\l}$ and let
$C=\set{\a:M_\a\cap N=N_\a}$.  If $\a\in C\cap S$ then $S\cap\a\in I_\a^+$
since otherwise there is $A\in N_\a\sse M_\a$ such that $S\sse A$ and
$\a\notin
A$, which is impossible.  [Note: we really only know that $S\cap\a\sse
A\cap\a$, but since $S,A\in M_\a$ we must have $M_\a\models
S\sse A$ and since $M_\a\prec H(\l^+)$, it is true that $S\sse A$.  We will
leave such remarks tacit in future.]

\medskip

For $i<k$ let $f_i(\a)=\set{\b<\a:f\{\b,\a\}=i}$.

Next, if $\s$ is a sequence of length $n$ of elements of $k$ we define an
ideal
$I(\a,\s)$ for $\a\in S_n$, where $S_{n+1}=\set{\a\in S_n: \a\hbox{\ is a
limit
point of }S_n}$.  The definition is by induction on $n=$ length($\s$).  Let
$I(\a,\langle\rangle)=I_\a$.  If $\s=\t^\frown\langle i\rangle$ then put 
$$X\in
I(\a,\s)\quad{\rm iff}\quad \set{\b<\a: X\cap f_i(\a)\cap\b\in I^+(\b,\t)}\in
I_\a.$$ Notice
that some of the $I(\a,\s)$ may not be proper.  As part of the proof of the
Erd\H os-Rado Theorem we essentially made the easy observation that there
exists $i$ such that $I(\a,\langle i\rangle)$ is proper.

\proclaim Lemma 3.3. We always have $I_\a\sse I(\a,\s)$.

\proof Once again we proceed by induction on length($\s$).  Suppose
$\s=\t^\frown\langle i\rangle$.  We may assume $I(\a,\s)$ is proper.  Let
$X\in
I_\a$.  then $\exists A\in N_\a$ $X\sse A\cap\a$ and $\a\notin A$.  We may as
well assume $X=A\cap\a$.  Since $\a$ is limit and $A\in N_\a$ we know
$\exists\g<\a$ $A\in N_\g$.  Suppose $\b<\a$, $\g<\b$, and $A\cap\b\in
I^+(\b,\t)$.  Then by inductive hypothesis $A\cap\b\in I_\b^+$ so $\b\in A$. 
Thus
$$\set{\b<\a:(A\cap\a)\cap f_i(\a)\cap\b\in
I^+(\b,\t)}\sse\set{\b<\a:A\cap\b\in
I^+(\b,\t)}\sse A\cap\a\in I_\a,$$
so $A\cap\a\in I(\a,\s)$ as desired.

\medskip

Now we extend an observation made in the proof of the Erd\H os-Rado Theorem.

\proclaim Lemma 3.4. If $X\sse\a$ and $X\in I^+(\a,\s)$ then $\forall j\in
{\rm
range}\,\s$ $\exists W\sse X$ $W$ is $j$-homogeneous and $|W|=\k$.

\proof The proof is by induction on length($\s$).  We may assume
$\s=\t^\frown\langle i\rangle$.  If $j\in{\rm range}(\t)$ then we may apply
the
inductive hypothesis to $X\cap\b\in I^+(\b,\t)$ for some $\b<\a$.  So suppose
$j=i$.  As in Lemma 2.2 let us argue that if $x\sse X$, $|x|<\k$ and
$x\cup\{\a\}$ is $i$-homogeneous then $\exists\x\in X$ $\x>\sup x$ and
$x\cup\{\x\}\cup\{\a\}$ is $i$-homogeneous.  Let
$A=\set{\b:x\cup\{\b\}\hbox{\ is
$i$-homogeneous}}$.  Then $A\in N_\a$ and $\a\in A$ so $a-A\in I_\a$.  By
Lemma~3.3 $\a-A\in I(\a,\s)$, so $A\cap X\in I^+(\a,\s)$.  Hence there is
$\b<\a$, $\b>\sup x$, such that $A\cap X\cap f_i(\a)\cap\b\in I^+(\b,\t)$. 
Now
just choose $\x\in A\cap X\cap f_i(\a)$ with $\x>\sup x$.

\medskip

Given an ordinal $\r$ and a finite sequence $\s$ of elements of $k$ we define
what it means for $x\sse\l$ to be $(\r,\s)$-{\it good}.  $x$ is
$(\r,\langle\rangle)$-good iff $x$ is a singleton.  If $\s=\t^\frown\langle
i\rangle$ then $x$ is $(\r,\s)$-good iff $x=\bigcup\set{x_\x:\x<\r}$ where
each
$x_\x$ is $(\r,\t)$-good, $\x<\h<\r$ implies $\sup x_\x<\inf x_\h$, and
whenever $\x<\h<\r$, $\g\in x_\x$ and $\d\in x_\h$ then $f\{\g,\d\}=i$.  Note
that any $(\r,\s)$-good set must have order type $\r^n$ where $n=$
length($\s$).

\proclaim Lemma 3.5. If $x$ is $(\r,\s)$-good then $\forall i\in{\rm
range}\,\s$ $\exists y\sse x$ $y$ is $i$-homogeneous and has order type $\r$.

\proof Easy by an induction on length($\s$).

\proclaim Lemma 3.6. If $X\in I^+(\a,\s)$ then $\forall\r<\k$ $\exists x\sse
X$
$x$ is $(\r,\s)$-good.

\proof Say $\s=\t^\frown\langle i\rangle$.  We construct
$x=\bigcup\set{x_\x:\x<\r}\sse X\cap f_i(\a)$.  By induction on $\x$ we
obtain
$(\r,\t)$-good $x_\x\sse X\cap f_i(\a)$.  Suppose $x_\h$ has been obtained
for $\h<\x$.  Since $\r<\k$ we know $x_\h\in N_\a$ for each $\h<\x$ and hence
$\seq{x_\h:\h<\x}\in N_\a$.  Thus $A=\set{\b:\forall\h<\x\ \forall\g\in x_\h\
f\{\g,\b\}=i}\in N_\a$.  Moreover, $\a\in A\in I_\a^*\sse I^*(\a,\s)$, so
$X\cap
A\in I^+(\a,\s)$.  Find $\b<\a$ such that $X\cap A\cap f_i(\a)\cap\b\in
I^+(\b,\t)$ and choose $(\r,\t)$-good $x_\x\sse X\cap A\cap f_i(\a)\cap\b$. 
We may assume $\b$ is large enough so that we may choose $x_\x$ with $\sup
x_\h<\min x_\x$ for all $\h<\x$.  The rest is easy.

\medskip

Recall that an {\it indecomposable\/} ordinal is an ordinal power of $\w$. 
An
indecomposable ordinal is characterized by the property that whenever
$\r=A\cup
B$, either $A$ or $B$ has order type $\r$.

\proclaim Lemma 3.7. If $\r$ is indecomposable, $x\sse\l$ is $(\r,\s)$-good,
$m<\w$ and $x=\bigcup\set{t_j:j<m}$ then $\exists j$ $\exists y\sse t_j$ $y$
is
$(\r,\s)$-good.

\proof Let $\s=\t^\frown\langle i\rangle$ and let $\seq{x_\x:\x<\r}$
witness $(\r,\s)$-goodness of $x$.  By inductive hypothesis for each $\x<\r$
there is $j(\x)<m$ and $y_\x\sse x_\x\cap t_{j(\x)}$ such that $y_\x$ is
$(\r,\t)$-good (this is trivial if $\t=\langle\rangle$).  Since $\r$ is
indecomposable there is $j$ such that $\set{\x:j(\x)=j}$ has order type $\r$. 
Let $y=\bigcup\set{y_\x:j(\x)=j}$.  Then $y\sse t_j$ and $y$ is
$(\r,\s)$-good.

\medskip

If $\s$ is a finite sequence of elements of $k$ let
$f_\s(\a)=\bigcup\set{f_i(\a): i\in{\rm range}\,\s}$.

Recall that we defined $S_{n+1}=\set{\a\in S_n:\a\hbox{\ is a limit point of
}S_n}$.  Thus $S_\w=\bigcap\set{S_n:n<\w}$ differs from $S_0$ by a
nonstationary set, and $I(\a,\s)$ is defined for all $\s$ whenever $\a\in
S_\w$.

If $T\sse S_\w$ is stationary and $\a\in T$, let $\Sigma(\a,T)=\{\,\s:\s$ is
one-to-one, length($\s$) $\ge 1$, and $T\cap\a\in I^+(\a,\s)\,\}$.  Note that
if $\a$ is a reflection point of $T$ then $\Sigma(\a,T)\ne0$.

\proclaim Lemma 3.8. There is stationary $S\sse S_\w$ and $\Sigma$ such that
for all stationary $T\sse S$ there is a closed unbounded set $C$ such that
$\forall\a\in T\cap C$ $\Sigma(\a,T)=\Sigma$.

\proof Note that if $\a\in T\sse S$ then $\Sigma(\a,T)\sse\Sigma(\a,S)$. 
Choose descending sequences
$$S_\w=T_0\supseteq T_1\supseteq T_2\supseteq\ldots\quad{\rm and}\quad
\Sigma_0\supseteq\Sigma_1\supseteq\Sigma_2\supseteq\ldots$$
such that $\forall n$ $\forall\a\in T_{n+1}$ $\Sigma(\a,T_n)=\Sigma_n$, the
$T_n$ are stationary and $\Sigma_{n+1}$ is a proper subset of $\Sigma_n$. 
Since each $\Sigma_n$ is finite this process must end after finitely many
steps. 
Let $S$ and $\Sigma$ be the final elements of each sequence.

\medskip

Clearly $\Sigma\ne0$.  Choose $\s\in\Sigma$ maximal with respect to
inclusion.

\proclaim Lemma 3.9. There are $\a\in S$ and stationary $T\sse S$ such that
$\forall\b\in T$ $\a-f_\s(\b)\in I(\a,\s)$.

\proof Suppose the Lemma is false.  Then for each $\a\in S$ there is closed
unbounded $C_\a$ such that $\forall\b\in C_\a\cap S$ $\a-f_\s(\b)\in
I^+(\a,\s)$.

Let $C=\set{\b:\forall\a\in\b\cap S\ \b\in C_\a}$, the diagonal intersection
of
the $C_\a$.  Then $C$ is closed unbounded.  Let $\a\in S$ be a reflection
point
of $C\cap S$ with $\s\in\Sigma(\a,S)$.  For all $\b\in C\cap S\cap\a$ there
must be $i\notin$ range($\s$) such that $\b\cap f_i(\a)\in I^+(\b,\s)$.  Fix
$i<k$ such that $\set{\b\in C\cap S\cap\a:\b\cap f_i(\a)\in I^+(\b,\s)}\in
I_\a^+$.  But now clearly if $\t=\s^\frown\langle i\rangle$ then $S\cap \a\in
I^+(\a,\t)$, contradicting maximality of $\s$.

\medskip

Now we are ready to put all this material together to complete the proof of
Theorem 3.1.

Fix an indecomposable ordinal $\r<\log\k$.  It will suffice to find a
homogeneous set of order type $\k+\r$.  Let $\a$ and $T$ be as in Lemma 3.9. 
Choose $\b>\a$ such that $\Sigma(\b,T)=\Sigma$.  Since $T\cap\b\in
I^+(\b,\s)$
there is $x\sse T$ that is $(\r,\s)$-good.  We may assume $\min x>\a$.  Also,
since $I(\a,\s)$ is $\k$-complete we know $F=\bigcup\set{\a-f_\s(\b):\b\in
x}\in
I(\a,\s)$.  For each $\g\in\a-F$ define $g_\g:x\to$ range($\s$) by
$g_\g(\b)=f\{\g,\b\}$.  Since $|x|<\log\k$ there is $g$ such that
$\set{\g:g_\g=g}\in I^+(\a,\s)$.  Also, by Lemma 3.7 $\exists i$ $g^{-1}(i)$
is
$(\r,\s)$-good.

But now by Lemma 3.4 $\exists W\sse\set{\g:g_\g=g}$ such that $W$ is
$i$-homogeneous, and by Lemma 3.5 we may find $y\sse g^{-1}(i)$ such that $y$
is $i$-homogeneous.  Also, if $\g\in W$ and $\b\in y$ than
$f\{\g,\b\}=g_\g(\b)=g(\b)=i$ so $W\cup y$ is $i$-homogeneous of order type
$\k+\r$, and the proof is complete.

\bigskip

\noindent{\bf 4. An unbalanced extension of the Erd\H os-Rado Theorem.}  It
is
natural to ask if Theorem 3.1 can be improved in the same way that Theorem
2.1
was improved to Theorem 2.3.  For example, is it true under CH that
$\w_2\to(\w_2,\w_1+2)^2$?  Hajnal showed in [Ha] that if
GCH holds and $\k$ is regular, then $\k^+\not\to(\k^+,\k+2)^2$, and Todor\v
cevi\'c has shown in unpublished work that this remains the case when $\k$ is
singular.  Thus the best we could hope to prove under CH would be
$\w_2\to(\a,\w_1+2)^2$ for all $\a<\w_2$.  In this section we will prove (a
generalization of) an improvement of part of this.

\proclaim Theorem 4.1. Suppose $\k$ is regular and $\l=(2^{<\k})^+$.  Then
$\l\to(\r,(\k+n)_k)^2$ for all $n,k<\w$, where $\r=\k^{\w+2}+1$.

Here $\k^{\w+2}$ represents ordinal exponentiation.  Recall that the
partition
relation means that if $f:[\l]^2\to k+1$, then either there is a
0-homogeneous
set of order type $\r$ or else there is an $i$-homogeneous set of type $\k+n$
for some $i>0$.

The rest of this section will be devoted to the proof of Theorem 4.1.  The
strategy of the proof is to derive Theorem 4.1 from the auxiliary assumption
$Q(\k)$, which asserts that $2^{<\k}=\k$ and in addition that
$$\forall\seq{f_\a:\a<\k^+}\sse{}^\k\k\ \exists g\in{}^\k\k\ \forall\a\
\exists
\h<\k\ \forall\x>\h\ f_\a(\x)<g(\x).\leqno(*)$$
Then, rather as in [BH], we observe that the assumption $Q(\k)$ is
unnecessary,
and therefore that Theorem 4.1 holds in ZFC.

Let us deal with this latter observation first.

Let $P_0$ be the natural $\k$-closed ordering for making $2^{<\k}=\k$.  Then
in
$V^{P_0}$ we have $\l=\k^+$.  Working in $V^{P_0}$ and using a standard
iterated forcing argument (as in [Ba]) we can force (*) to be true via a
partial
ordering $P_1$ that is $\k$-closed and has the $\l$-chain condition.  Let
$P=P_0*P_1$.  Then $P$ is $\k$-closed and in $V^P$ $\l=\k^+$ and $Q(\k)$
holds.  Note that in $V^P$ we will have $2^\k>\k^+$ (since for one thing that
is implied by (*)).

Assuming we have proved Theorem 4.1 under the assumption $Q(\k)$ we may
assume
it holds in $V^P$.  Assume $f:[\l]^2\to k+1$, where $f\in V$.  Then in $V^P$
there is $A\sse\l$ such that either (a) $A$ is 0-homogeneous of type $\r$ or
(b) $A$ is $i$-homogeneous of type $\k+n$, where $i>0$.  Suppose (a) holds. 
Note that $\k^{\w+2}+1$ is the same whether computed in $V$ or in $V^P$.  Let
$h:\k\to\r$ be a bijection with $h\in V$.  In $V^P$ fix an order-isomorphism
$j:\r\to A$.  Now, working in $V$, find a decreasing sequence
$\seq{p_\x:\x<\k}$ of elements of $P$ and a sequence $\seq{\a_\x:\x<\k}$ of
elements of $\l$ such that $\forall\x$ $p_\x\forces j(h(\x))=\a_\x$.  This is
easy to do by induction on $\x$, using the fact that $P$ is $\k$-closed.  But
now it is clear that $\set{\a_\x:\x<\k}\in V$ has order type $\r$ and is
0-homogeneous.  Case (b) may be handled the same way.

From now on, assume $Q(\k)$ holds (so $\l=\k^+$).

We begin with a couple of observations about order types.  The only use of
$Q(\k)$ will be in Corollary 4.3 below.

\proclaim Lemma 4.2. Let $n<\w$, $n\ge1$, $\a<\k^+$, and assume $\a$ has a
cofinal subset of type $\k^n$.  Then for every $f:\k\to\k$ there is a set
$A(f,n,\a)$ of order type $\k^n$ cofinal in $\a$ such that
\item{(i)} if $A$ is cofinal in $\a$ of type $\k^n$ then $\exists f$ $A\sse
A(f,n,\a)$
\item{(ii)} $\forall f_1, f_2$ if $\forall\x<\k$ $f_1(\x)\le f_2(\x)$ then
$A(f_1,n,\a)\sse A(f_2,n,\a)$.

\proof The proof is by induction on $n$.  Suppose $n=1$.  Let
$\seq{\a_\x:\x<\k}$ be an increasing, continuous sequence cofinal in $\a$
with
$\a_0=0$.  Let $A_\x=\a_{\x+1}-\a_\x$ and let $h_\x:\k\to A_\x$ be onto. 
Define $A(f,n,\a)=\bigcup\set{h_\x``f(\x):\x<\k}$.  It is easy to see that
this
works.

Now suppose $n>1$. Let $\seq{\a_\x:\x<\k}$ be as above.  If
$\a_\x<\b\le\a_{\x+1}$ then for $1\le i<n$ let
$A'(f,i,\b)=A(f,i,\b)\cap(\a_{\a+1}-\a_\x)$ if $A(f,i,\b)$ is defined; let
$A'(f,i,\b)=\{\b\}$ otherwise.  Let $A(f,0,\b)=\{\b\}$.  Now let $f:\k\to\k$. 
We must define $A(f,n,\a)$.  We may regard $f$ as being defined on
$\k\times\k$, i.e., $f:\k\times\k\to\k$, and we define $f_\x$ by
$f_\x(\h)=f(\x,\h)$ for $\x,\h<\k$.  Let $h:\k\to\a$ be a bijection.  Now let
$$A(f,n,\a)=A(f_0,1,\a)\cup\bigcup\set{A'(f_\x,i,h(\x)):h(\x)\in A(f_0,1,\a),
i<n}.$$

Since $A(f_0,1,\a)$ has type $\k$, it is clear that
$\bigcup\set{A'(f_\x,i,h(\x)): h(\x)\in A(f_0,1,\a)\cap(\a_{\h+1}-\a_\h)}$
has
type $<\k^n$ for fixed $\h$.  Thus $A(f_0,1,\a)$ has type at most $\k^n$.

It is clear that condition (ii) is satisfied.  Let us check (i).  Fix
$A\sse\a$ cofinal of type $\k^n$.  For $\h<\k$ let
$A_\h=A\cap(\a_{\h+1}-\a_\h)$.  Then $A_\h$ has type $<\k^n$ so there is a
family $F_\h$ such that $|F_\h|<\k$, $A_\h=\bigcup F_\h$, every $B\in F_\h$
has type $\k^i$ for some $i<n$ (possibly $i=0$), and if $B,C\in F_\h$ and
$B\ne C$ then $\sup B<\inf C$ or $\sup C<\inf B$.  Let $H_\h=\set{\sup B:B\in
F_\h}$.  Now choose $f:\k\times\k\to\k$ so that $\bigcup\set{H_\h:\h<\k}\sse
A(f_0,1,\a)$ and whenever $B\in F_\h$, $B$ has type $\k^i$ for $i>0$, and
$\sup B=\b=h(\x)$, then $B\sse A(f_\x,i,\b)$.  But now it is clear that
$A\sse
A(f,n,\a)$.

There may be special circumstances in which the $A(f,n,\a)$ as defined above
have order type $<\k^n$.  This defect may be remedied by choosing a fixed set
of order type $\k^n$ cofinal in $\a$ and adjoining its elements to every
$A(f,n,\a)$.

\proclaim Corollary 4.3. ($Q(\k)$)  Let $\a<\k^+$ and let $A_\g$, $\g<\k^+$,
be
a sequence of sets of order type $\k^n$ cofinal in $\a$, where $1\le n\le\w$. 
Then we may write $\k^+=\bigcup\set{X_\x:\x<\k}$ where
$\bigcup\set{A_\g:\g\in
X_\x}$ has order type $\k^n$ for each $\x$.

\proof First suppose $n<\w$.  For each $\g$ choose $f_\g:\k\to\k$ such that
$A_\g\sse A(f_\g,n,\a)$.  Use $Q(\k)$ to find $f:\k\to\k$ such that
$\forall\g$
$\exists\x_\g$ $\forall\h>\x_\g$ $f_\g(\h)<f(\h)$.  Thus there is
$s_\g:\x_\g\to\k$ such that if $g_\g=s_\g\cup f\rest(\k-\x_\g)$ then
$\forall\h$
$f_\g(\h)<g_\g(\h)$.  By $\k^{<\k}=\k$ there are only $\k$ functions $s_\g$
so
we may write $\k^+=\bigcup\set{X_\x:\x<\k}$ where if $\g,\d\in X_\x$ then
$s_\g=s_\d$, so $g_\g=g_\d$.  But then
$$\bigcup\set{A_\g:\g\in X_\x}\sse\bigcup\set{
A(f_\g,n,\a):\g\in X_\x}\sse A(g_{\g_0},n,\a)$$
by Lemma 4.2(ii) where $\g_0$ is any member of $X_\x$.

Now suppose $n=\w$.  For each $\g<\k^+$ let $g_\g(0)=0$ and if $1\le m<\w$
let
$g_\g(m)$ be the supremum of the first $\k^m$ elements of $A_\g$ and let
$A^m_\g=A_\g\cap[g_\g(m-1),g_\g(m))$.  Since $\k^\w=\k$ we may write
$\k^+=\bigcup\set{Y_\x:\x<\k}$ where for each $\g,\d\in Y_\x$ we have
$g_\g=g_\d$.  And now by the earlier case, for each $m$ we may write
$Y_\x=\bigcup\set{Y^m_{\x\h}:\h<\k}$, where $\bigcup\set{A^m_\g:\g\in
Y^m_{\x\h}}$ has order type $\k^m$ for each $\x$ and $\h$.  For each $\g$ let
$h_\g(m)=(\x,\h)$ if $\g\in Y^m_{\x\h}$.  Finally put
$\k^+=\bigcup\set{X_\x:\x<\k}$ where for $\g,\d\in X_\x$ we have $h_\g=h_\d$. 
It is easy to check that this works.

\medskip
This is the only use of $Q(\k)$ that we need.

For the main proof let us fix a partition function $f:[\l]^2\to k+1$, and
assume that there is no $i$-homogeneous set of type $\k+n$, where $i>0$.  We
must find a 0-homogeneous set of type $\r$.

If $H(\l^+)$ is the collection of sets hereditarily of cardinality $\le\l$,
then $f\in H(\l^+)$.  Let $S=\set{N:N\prec H(\l^+),\ |N|=\k,\ N\cap\l\in\l,\
f\in N,\ {\rm and}\ [N]^{<\k}\sse N}$.  Then $S$ is a stationary subset of
$[H(\l^+)]^\k$.

Define $\pi:S\to\l$ by $\pi(N)=N\cap\l$.  If $T\sse S$ is stationary then let
$J_T=\set{X\sse\l:\pi^{-1}(X)\cap T\ \hbox{is nonstationary}}$.  It is easy
to
see that $J_T$ is a normal ideal on $\l$.

\proclaim Claim 4.4. Let $\x<\w+2$.  Then
\item{(a)} If $T\sse S$ is stationary, then there is a set $C$ closed and
unbounded in $[H(\l^+)]^\k$ such that for all $N\in C\cap T$, if $A\in N$ and
$\pi(N)\in A$ then there is $x\in N$ such that $x\sse A\cap\pi(T)$, $x$ is of
order type $\k^\x$, and $x\cup\{\pi(N)\}$ is 0-homogeneous.
\item{(b)} If $T\sse S$ is stationary then there is stationary $U\sse T$ and
$x\sse\pi(T)$ such that $x$ is 0-homogeneous of type $\k^\x$ and
$\forall\a\in
x$ $\forall\b\in\pi(U)$ $f\{\a,\b\}=0$.

Let us observe that Claim 4.4(a) for $\x=\w+1$ will complete the proof.  Fix
$N\in C\cap S$.  By induction on $\g<\k$ we define $A_\g\in N$ with
$\pi(N)\in
A_\g$ and $x_\g\sse A_\g$ such that $x_\g$ has order type $\k^{\w+1}$ and
$x_\g\cup\{\pi(n)\}$ is 0-homogeneous.  Let $A_0$ be arbitrary and choose
$x_0$
as in 4.4(a).  Suppose $A_\g,\ x_\g$ have been determined for $\g<\d$.  Since
$N\in S$ we know $[N]^{<\k}\sse N$ so $\seq{x_\g:\g<\d}\in N$  Thus
$$A_\d=\left\{\a:\bigcup\set{x_\g:\g<\d}\cup\{\a\}\ \hbox{is 0-homogeneous
and }
\sup\bigcup\set{x_\g:\g<\d}<\a\right\} \in N$$
and $\pi(N)\in A_\d$.  Choose $x_\d\sse A_\d$ as in 4.4(a).  This completes
the
construction.  And now $\bigcup\set{x_\g:\g<\k}\cup\{\pi(N)\}$ is
0-homogeneous
and has order type $\k^{\w+1}\cdot\k+1=\k^{\w+2}+1=\r$, as desired.

Thus we may devote the rest of the section to Claim 4.4.  The proof is by
induction on $\x$.

Let us begin the induction by showing that (b) holds for $\x=0$.  If there is
a
$\a\in\pi(T)$ such that $X_\a=\set{\b:f\{\a,\b\}>0}\in J_T^+$, then let
$x=\{\a\}$ and $U=T\cap\pi^{-1}X_\a$.  If not, then the diagonal intersection
$X=\set{\b:\forall \a<\b\ \b\notin X_\a}$ belongs to $J_T^+$, hence has
cardinality $\l$.  But now by Theorem 3.1 there is an $i$-homogeneous set of
type $\k+n$ for some $i>0$, contrary to hypothesis.

Next we show (b) implies (a).  If (a) is false, then there is stationary
$S'\sse S$ such that the assertion in (a) fails for all $N\in S'$.  By
normality of the nonstationary ideal on $[H(\l^+)]^\k$ there is a single set
$A$ and stationary $T\sse S'$ such that $A$ is a counterexample to (a) for
all
$N\in T$.  Since $\pi(N)\in A$ for $N\in T$ we have $\pi(T)\sse A$.  Let $x$
and $U$ be as in (b).  Then $x\sse\pi(T)\sse A\cap\pi(S')$ and $x$ is of type
$\k^\x$.  Since $x\in H(\l^+)$ there must be $N\in U$ with $x\in N$. 
Moreover
$x\cup\{\pi(N)\}$ is 0-homogeneous, so $A$ is not in fact a counterexample
for
$N$, contradiction.

We now concentrate on showing that (a) for $\x$ implies (b) for $\x+1$,
provided $\x\le\w$.

Let $T$ be as in (b), let $C$ be as in (a), and let $\seq{N_\a:\a<\l}$ be an
increasing continuous sequence of elements of $C$ with $[N_\a]^{<\k}\sse
N_{\a+1}$.  Let $N_\l=\bigcup\set{N_\a:\a<\l}$.  We may assume that for all
$N\in T$, if $\pi(N)=\a$ then $N\cap N_\l=N_\a$.  (If $(N\cap
N_\l)-N_{\pi(N)}\ne 0$ for stationarily many $N$, then for stationarily many
$N$ there is a common element of $(N\cap N_\l)-N_{\pi(N)}$, and this is
impossible.)  We also assume that if $N\in T$ and $N\cap N_\l=N_\a$, then
$N_\a\cap\l=\a$.

Let $\Sigma_0$ be the set of all one-to-one sequences $\s$ of elements of
$k+1$
such that $\s(0)=0$.  For each $\a\in\pi(C\cap T)$ and each $\s\in\Sigma_0$
we
will define an ideal $I(\a,\s)$ on $\a$.  Each $I(\a,\s)$ will contain all
bounded subsets of $\a$ but will not necessarily be closed under countable
unions.

First suppose $\s=\langle0\rangle$.  Choose $N\in T\cap C$ with $\pi(N)=\a$. 
Let $\seq{A_\g:\g<\k}$ enumerate $\set{A\in N:\a\in A}$.  We define
$\seq{x_\g:\g<\k}$ by induction on $\g$ so that $x_\g\in N$, $x_\g$ has type
$\k^\x$ and $x_\g\cup\{\a\}$ is
0-homogeneous.  Suppose $x_\g$, $\g<\d$, has been obtained.  Let
$$A=\bigcap\set{A_\g:\g<\d}\cap\Big\{\,\b:\b>\sup\set{x_\g:\g<\d}\ \hbox{and
}\forall\g<\d\ \forall\b'\in x_\g\ f\{\b,\b'\}=0\,\Big\}.$$
Note that $A\in N$ since it is definable from $\seq{A_\g:\g<\d}$ and
$\seq{x_\g:\g<\d}$, both of which are in $N$, and clearly $\a\in A$.  By
4.4(a)
choose $x_\d\in N$, 0-homogeneous of type $\k^\x$, such that $x_\d\sse
A\cap\pi(T)$.  Thus $\bar x(\a)=\bigcup\set{x_\g:\g<\k}$ is 0-homogenous of
type $\k^{\x+1}$, and is cofinal in $\a$ since for any $\b<\a$ $\l-\b=A_\g$
for
some $\g$.  Let $I(\a,\langle0\rangle)=\set{X\sse \a:X\cap\bar x(\a)\
\hbox{has
type }<\k^{\x+1}}$.

Now suppose $\t\in\Sigma_0$ and $\s=\t^\frown\langle i\rangle$.  Let $I_\a$
be
the ideal on $\a$ generated by $\set{A\cap\a:A\in N_\a, \a\notin A}$, as in
the
previous section.  Recall that $f_i(\a)=\set{\b<\a:f\{\b,\a\}=i}$.  For
$X\sse
\a$, put $X\in I(\a,\s)$ iff $\set{\b<\a:X\cap f_i(\a)\cap\b\in
I^+(\b,\t)}\in
I_\a$.

Next we prove analogues of Lemmas 3.3, 3.4, and 3.9.

\proclaim Lemma 4.5.  $I_\a\sse I(\a,\s)$ for $\a\in\pi(C\cap T)$,
$\s\in\Sigma_0$.

\proof Let $X\in I^+(\a,\s)$.  If $\s=\langle0\rangle$, then by construction
$X\cap A\ne 0$ for all $A\in N_\a$ with $\a\in A$.  Thus $X\in I^+_\a$. 
Suppose $\s=\t^\frown\langle i\rangle$, where $\t\in\Sigma_0$. Let $A\in
N_\a$,
$\a\in A$.  Since $\a$ is limit $\exists\g<\a$ $A\in N_\g$.  Choose $\b$ such
that $\g<\b<\a$, $\b\in A$, and $X\cap\b\in I^+(\b,\t)$.  By inductive
hypothesis
$X\cap\b\in I_\b^+$, so $X\cap\b\cap A\ne0$.  Thus $X\cap A\ne0$ so $X\in
I_\a^+$.

\proclaim Lemma 4.6.  If $X\in I^+(\a,\s)$, then for all $j\in \range\s$, if
$j>0$ then $X$ contains a $j$-homogeneous set of type $\k$.

\proof We proceed by induction.  Say $\s=\t^\frown\langle i\rangle$ where
$\t\in\Sigma_0$.  The induction will handle the case
$j\in\range\t$.  It will suffice to show, as in Section 3, that if
$x\cup\{\a\}$ is $i$-homogeneous and $|x|<\k$ then $\exists\b\in X-x$
$x\cup\{\b,\a\}$ is $i$-homogeneous.  Let $A=\{\,\g:x\cup\{\g\}$ is
$i$-homogeneous$\,\}$.  Then
$A\in N_\a$ and $\a\in A$.  Since $X\in I^+(\a,\s)$ there is $\g\in A$ such
that $X\cap f_i(\a)\cap\g\in I^+(\g,\t)$.  We may choose $\g$ large enough
that
$A\in N_\g$ and $\g>\sup x$.  Thus $X\cap f_i(\a)\cap\g\cap A\ne0$ and we may
choose $\b\in X\cap f_i(\a)\cap\g\cap A$, $\b>\sup x$.

For $\b<\l$ and $\s\in\Sigma_0$, let
$f_\s(\b)=\bigcup\set{f_i(\b):i\in\range\s}$.

\proclaim Lemma 4.7.  There are $\s\in\Sigma_0$, $\a\in\pi(C\cap T)$, and
stationary $T'\sse T$ such that $\a\in I^+(\a,\s)$ (i.e., $I(\a,\s)$ is
proper)
and $\forall\b\in\pi(T')$ $\a-f_\s(\b)\in I(\a,\s)$.

\proof For $\a\in\pi(C\cap T)$ we know $\pi(T)\in
I^+(\a,\langle0\rangle)$.  Choose $\s_\a\in\Sigma_0$ maximal with respect to
inclusion such that $\pi(T)\in I^+(\a,\s)$, and let $U\sse C\cap T$ be
stationary such that $\forall\a\in\pi(U)$ $\s_\a=\s$.

Suppose the Lemma is false.  Then for each $\a\in\pi(U)$ there is $C_\a\in
J_U^*$ such that $\forall\b\in C_\a$ $\a-f_\s(\b)\in I^+(\a,\s)$.  (Recall
that
$J_U=\set{X\sse\l:\pi^{-1}(X)\cap U\hbox{\ is stationary}}$.)  Since $J_U$ is
normal, the set $C_\l=\set{\b:\forall\a<\b\ \b\in C_\a}\in J_U^*$ also.  By
Lemma 3.2 there is closed unbounded $D\sse\l$ such that $\forall\a\in
D\cap\pi(U)$\ \ $\pi(U)\cap\a\in I_\a^+$.  Fix $\a\in D\cap C_\l\cap\pi(U)$.

If $\b<\a$ and $\b\in\pi(U)$ then $\a\in C_\b$ so there is $i$,
$i\notin\range\s$, such that $f_i(\a)\cap\b\in I^+(\b,\s)$.  Since $k$ is
finite we may fix $i$ such that $\set{\b<\a:\b\in\pi(U)\ {\rm and}\
f_i(\a)\cap\b\in I^+(\b,\s)}\in I_\a^+$.  But this means that
$I(\a,\s^\frown\langle i\rangle)$ is proper contrary to the assumption about
the maximality of $\s$.

One more lemma will allow us to complete this part of the proof of Claim 4.4.

\proclaim Lemma 4.8.  Let $\a<\l$ and $\s\in\Sigma_0$ be arbitrary.  If $X\in
I^+(\a,\s)$, $Y\in I(\a,\s)$ and $Y\sse X$, then $\exists \b\le\a$ $\bar
x(\b)\cap X$ has type $\k^{\x+1}$ and $Y\cap\bar x(\b)$ has type
$<\k^{\x+1}$.

\proof If $\s=\langle0\rangle$ this is clear.  Just take $\b=\a$.  Suppose
$\s=\t^\frown\langle i\rangle$ for $\t\in\Sigma_0$.  Then we can find $\b<\a$
such that $X\cap f_i(\a)\cap\b\in I^+(\b,\t)$ and $Y\cap\b\in I(\b,\t)$, so
we
are done by inductive hypothesis.

\medskip
Finally, let $\a$, $\s$ and $T'$ be as in Lemma 4.7.  For $\b\in\pi(T')$ let
$g(\b)=\bigcup\set{f_i(\b)\cap\a: i\in\range\s, i>0}$.  Let $F\sse \pi(T')$
be
maximal such that $\set{\cap\set{g(\b):\b\in x}:x\in[F]^{<\w}}\sse
I^+(\a,\s)$.  (It is possible that $F=0$.)

First, suppose $|F|=\l$.  Let $G$ be an ultrafilter on $\a$ such that
$\set{g(\b):\b\in F}\sse G\sse I^+(\a,\s)$.  Thus for each $\b\in F$ there is
$i_\b$ such that $f_{i_\b}(\b)\cap\a\in G$.  Find $i$ so that $F'=\set{\b\in
F:
i_\b=i}$ has cardinality $\l$.  We may assume we are proving Theorem 4.1 by
induction on $k$, so by $\l\to(\l,\k)^2$, which follows from Theorem 2.3, we
may assume there is $x\sse F'$ $i$-homogeneous of cardinality $n$.  But now
the
set $\bigcap\set{f_i(\b)\cap\a:\b\in x}$ belongs to $I^+(\a,\s)$ so by Lemma
4.6 must contain an $i$-homogeneous set $H$ of type $\k$.  And now $H\cup x$
is
$i$-homogeneous of type $\k+n$, contrary to our hypothesis about the
partition
function $f$.  Thus $|F|<\l$.

From the maximality of $F$ it follows that if $N\in T'$ and $\pi(N)\notin F$
then there is a finite set $y_N\sse F$ such that $\bigcap\set{g(\b):\b\in
y_N}
- f_0(\pi(N))\in I(\a,\s)$.  Since $|F|<\l$ there is $y$ such that
$T''=\set{N\in
T':y_N=y}$ is stationary.  For $N\in T''$ let $a_N=\bigcap\set{g(\b):\b\in
y}-f_0(\pi(N))$.  By Lemma 4.8 there is $\b_N\le\a$ such that $\bar
x(\b_N)\cap\bigcap\set{g(\b):\b\in x}$ has type $\k^{\x+1}$ and $\bar
x(\b_N)\cap a_N$ has type $<\k^{\x+1}$.  There is stationary $T'''\sse T''$
and
$\b$ such that $T'''=\set{N\in T'':\b_N=\b}$.

By thinning out $T'''$ we may assume without loss of generality that each
$a_N$ for $N\in T'''$ is such that $a_N\sse(\bar x(\b)\cap\z)\cup b_N$ where
$\z<\sup \bar x(\b)$ and $b_N$ has type $\k^\h$ for some $\h\le\x$, and that
$\h$ and $\z$ are the same for all $N$.  Now it follows from Corollary 4.3
that
$T'''=\bigcup\set{T_\d:\d<\k}$ where for each $\d$, $\bigcup\set{a_N:N\in
T_\d}$
has type $<\k^{\x+1}$.  Pick such $\d$ for which $T_\d$ is stationary.  Let
$U=T_\d$ and let $x=\bar x(\b)\cap\bigcap\set{g(\b):\b\in
y}-\bigcup\set{a_N:N\in T_\d}$.  Then $x$ is 0-homogeneous of type
$\k^{\x+1}$
and $\forall\a\in x$ $\forall\b\in\pi(U)$ $f\{\a,\b\}=0$.

This completes the proof of Claim 4.4(b) for $\x+1$ from Claim 4.4(a) for
$\x$,
provided $\x\le\w$.

To complete the proof of Claim 4.4, and hence Theorem 4.1, let us observe
that
4.4(a) for $\x=\w$ follows from 4.4(a) for $\x<\w$.  Let $C_n$ be a closed
unbounded set as in 4.4(a) for $\x=n$, and let $C=\bigcap\set{C_n:n<\w}$. 
Let
$N\in C\cap T$, $A\in N$, and $\pi(N)\in A$.  By induction on $n$ choose
$x_n\in N$ as follows.  Let $x_0$ be as in 4.4(a) for $\x=0$.  Given $x_m$
for
$m<n$, let
$$A_n=\set{\b:\bigcup\set{x_m:m<n}\cup\{\b\}\ \hbox{is
0-homogeneous}}.$$
Then $A_n\in N$ and $\pi(N)\in A$, so we may find $x_n\in N$
as in 4.4(a) with $A$ replaced by $A_n$.  But now $x=\bigcup\set{x_n:n<\w}\in
N$ and satisfies 4.4(a) with $\x=\w$.

\bigskip

\noindent{\bf 5. An extension of the Erd\H os-Rado Theorem to ordinals.}  In
this section we prove a theorem implying that if CH holds then for any $n<\w$
there is $m<\w$ such that $\w_2\cdot m\to(\w_1\cdot n)^2_2$.  This is an
improvement of a result of Shelah [Sh2] that implies under similar
circumstances
$\w_2^m\to(\w_1\cdot n)^2_2$.

Let $A$ be a set and let $k<\w$.  A {\it set-mapping\/} on $A$ of order $k$
is
a function $p:A\to[A]^{<k}$ such that $a\notin p(a)$ for all $a\in A$.  Given
such a mapping $p$, a set $F\sse A$ is said to be {\it free\/} for $p$ if
$\forall a,b\in F$ $a\notin p(b)$.  It is well-known that if $n,k<\w$ then
there is $g(n,k)<\w$ such that any set-mapping on $g(n,k)$ of order $k$ has a
free set of cardinality $n$.  (This is easy to see using elementary Ramsey
theory, for example.)

Let $n,k<\w$.

We define $f(n,k)$ by induction on $k$.  Let $f(n,0)=1$, and
let $f(n,k+1)=g(n,f(n,k))$.

\proclaim Theorem 5.1. Let $\k$ be regular and let $\l=(2^{<\k})^+$.  Then
for
all $n,k<\w$ we have
$$\l\cdot f(n,k)\to(\k\cdot n,\k\cdot(k+1))^2.$$

\proof The proof is by induction on $k$.  If $k=0$ then $f(n,k)=1$ so we must
show $\l\to(\k\cdot n,\k)^2$, and this follows easily from Theorem 2.3, which
implies $\l\to(\l,\k)^2$.

So we may assume $k>0$.  Let $m=f(n,k)$, and suppose $h:[m\times\l]^2\to2$. 
Note that with the lexicographic order, $m\times\l$ has order type $\l\cdot
m$.

Choose $N_0\prec N_1\prec \cdots\prec N_{m-1}\prec H(\l^+)$ such that $N_i\in
N_{i+1}$, $|N_i|=2^{<\k}$, $N_i\cap\l\in\l$, and $[N_i]^{<\k}\sse N_i$ for
all
$i$.  For each $i$ let $\b_i=N_i\cap\l$.

Let $A_i=\{i\}\times\l$.  Without loss of generality we may assume that for
each $i$ there is $n_i$ such that $A_i$ contains no 1-homogeneous set of type
$\k\cdot(n_i+1)$ but for all $B\sse A_i$, if $|B|=\l$ then $B$ contains a
1-homogeneous set of type $\k\cdot n_i$.  Of course if $n_i\ge k+1$ we are
done, so we may assume $n_i\le k$ for each $i$.

By induction on $\x<\k$ we will define $a_{i\x}\in A_i\cap N_0$ for each
$i<m$,
and we will obtain a set-mapping $g_\x$ on $m$.

Fix $i$, and let $X(i,\x)=\set{a\in A_i:\forall\h<\x\ \forall j<m\
f\{a,a_{j\h}\}=f\{(i,\b_i),a_{j\h}\}}$.  Since each $a_{j\h}\in N_0$ and
$[N_0]^{<\k}\sse N_0$, we have $X(i,\x)\in N_0$ (note that the function
assinging $a_{j\h}$ to $f\{(i,\b_i),a_{j\h}\}$ must be in $N_0$).  Since
$(i,\b_i)\in X(i,\x)$ it is clear that $|X(i,\x)|=\l$.

Now by $\l\to(\l,\k)^2$ we may assume $X(i,\x)$ contains a 1-homogeneous set
$Y$ of order type $\k$ (so $n_i\ge1$ also), and we may assume $Y\in N_0$
since
$N_0\prec H(\l^+)$ and the statement that such a $Y$ exists is true in
$H(\l^+)$.  If $Y\in N_0$, then in particular $Y\sse N_0$ as well.

We define $Y=Y_0\supseteq Y_1\supseteq\cdots\supseteq Y_m$ with $|Y_j|=\k$,
as
follows.  Given $Y_j$, let $Y_{j+1}=\set{y\in Y_j:f\{(j,\b_j),y\}=0}$ if this
set
has cardinality $\k$; otherwise set $Y_{j+1}=Y_j$ and put $j\in g_\x(i)$. 
Finally, choose $a_{i\x}\in Y_m$.  Note that $a_{i\x}\in N_0$.

For the following lemmas let us assume for convenience that Theorem 5.1 is
false
for~$k$.

\proclaim Lemma 5.2. $i\notin g_\x(i)$.

\proof Suppose $i\in g_\x(i)$.  Let $C=\set{y\in Y_i:f\{(i,\b_i),y\}=0}$. 
Since
$|C|<\k$ we have $C\in[N_0]^{<\k}\sse N_0$.  Now $(i,\b_i)$ belongs to
$Z=\set{b\in A_i:\forall y\in Y_i-C\ f\{b,y\}=1}$.  But $Y_i\in N_i$ since it
is defined using $Y$ and the $\b_j$ for $j<i$ and all these belong to $N_i$. 
Thus $Z\in N_i$ and since $(i,\b_i)\in Z$ we must have $|Z|=\l$.  By our
hypothesis on $n_i$ there is $B\sse Z$, 1-homogeneous of type $\k\cdot n_i$,
and
we may assume $\min B>\max Y_i$.  But now $(Y_i-C)\cup B$ is 1-homogeneous of
type $\k\cdot(n_i+1)$, contradiction.  Hence $i\notin g_\x(i)$.

\proclaim Lemma 5.3. $|g_\x(i)|<f(n,k-1)$.

\proof Again, suppose otherwise.  For each $j\in g_\x(i)$, choose a set
$C_j\in[Y_j]^{<\k}$ such that $(j,\b_j)\in Z_j=\set{b\in A_j:\forall y\in
Y_j-C_j\ f\{b,y\}=1}$.  Then $Z_j\in N_j$, and $|Z_j|=\l$.

By the inductive hypothesis for Theorem 5.1 we know $\l\cdot
f(n,k-1)\to(\k\cdot n,\k\cdot k)^2$.  If we apply this to
$\bigcup\set{Z_j:j\in
g_\x(i)}$, which has order type $\ge \l\cdot f(n,k-1)$, either this set
contains a 0-homogeneous set of type $\k\cdot n$ (in which case we have
established Theorem 5.1) or else it contains a 1-homogeneous set $B$ of type
$\k\cdot k$.  But then $B\cup Y_m-\bigcup\set{C_j:j\in g_\x(i)}$ is
1-homogeneous of type $\k\cdot(k+1)$ and again we are done.  Thus
$|g_\x(i)|<f(n,k-1)$.

Since there are only finitely many possibilities for $g_\x$ we may assume
that
they are all the same.  Say $g_\x=g$ for all $\x<\k$.  By the definition of
$m=f(n,k)$ and Lemmas 5.2 and 5.3 we know there is a free set $F$ for $g$ of
cardinality $n$.

\proclaim Lemma 5.4. If $\x\ne\h$ and $j_1,j_2\in F$ then
$f\{a_{j_1\x},a_{j_2\h}\}=0$.

\proof Assume $\x<\h$, $i=j_1$, $j=j_2$.  Since $j\notin g(i)$ we must have
had
$a_{i\x}\in Y_m\sse Y_{j+1}=\set{y\in Y_j:f\{(j,\b_j),y\}=0}$ during the
choice
of $a_{i\x}$.  Hence $f\{(j,\b_j),a_{i\x}\}=0$, so necessarily
$f\{a_{j\h},a_{i\x}\}=0$ also.
\endproof

And now if $F=\{i_0,\ldots,i_{n-1}\}$ and $X_0,\ldots,X_{n-1}$ are disjoint
subsets of $\k$, each of cardinality $\k$, then $\set{a_{i_l\x}:l<n\ {\rm
and}\
\x\in X_l}$ is 0-homogeneous of type $\k\cdot n$.  This completes the proof
of
Theorem 5.1.

\proclaim Corollary 5.5. $\l\cdot\w\to(\k\cdot n)^2_2$ for all $n<\w$.

It is natural to ask whether the subscript in Theorem 5.1 can be improved. 
By
Theorem 3.1 we know that if $\r<\log \k$ then $\l\to(\k+\r)^2_3$.  If
$\log\k<\k$ then a $\k$-closed forcing construction yields the consistency of
$\l\not\to(\k+\log\k)^2_2$ as in [Pr] (and we may arrange $2^\k>\l$ as well),
and
Corollary 5.5 shows $\l^+\to(\k\cdot n)^2_2$ for all $n<\w$.  What about
$\l^+\to(\k+\log\k)^2_3\,$?  No ordinal below $\l^+$ will work in view of the
following.  Throughout, we assume $\k$ is regular and $\l=(2^{<\k})^+$.

\proclaim Proposition 5.6. If $\l\not\to(\k+\log\k)^2_2$ then
$\forall\a<\l^+$
$\a\not\to(\k+\log\k,\k+\log\k,\w)^2$.

\proof Suppose $f:[\l]^2\to 2$ witnesses $\l\not\to(\k+\log\k)^2_2$.  Define
$g:[\a]^2\to3$ as follows.  Let $\pi:\a\to\l$ be a bijection (we assume
$\l\le\a<\l^+$).  If $\x<\h<\a$ then let $g\{\x,\h\}=2$ if $\pi(\x)>\pi(\h)$
and let $g\{\x,\h\}=f\{\pi(\x),\pi(\h)\}$ otherwise.  It is easy to see that
this works.

\proclaim Proposition 5.7. $\l^+\to(\k\cdot n,\k\cdot n,\k+1)^2$ for all
$n<\w$.

\proof A proof exactly similar to Theorem 2.3 shows that
$\l^+\to(\l^+,\k+1)^2$.  Thus if $f:[\l^+]^2\to3$, either there is a
2-homogeneous set of type $\k+1$ or else there is a set $X$ of cardinality
$\l^+$ such that $f``[X]^2\sse\{0,1\}$.  And now in the latter case Theorem
5.1
may be applied to find a subset of $X$ homogeneous of type $\k\cdot n$.
\endproof

\noindent{\bf Question.} Does $\l^+\to(\k+\log\k,\k+\log\k,\k+2)^2$?
\medskip
Assuming CH, the simplest nontrivial case is
$\w_3\to(\w_1+\w,\ \w_1+\w,\ \w_1+2)^2$ (and $2^{\w_1}\ge\w_3$).

\vfill
\eject

\centerline{References}
\bigskip
\parindent=0pt
\parskip=\medskipamount
\def\pp#1]#2\par{#1]\quad#2\par}
\pp[Ba]  J. Baumgartner, {\it Iterated forcing},
{\bf Surveys in Set Theory}, A. R. D. Mathias, ed., {\bf London Math. Soc.
Lecture Note Series 87} (1983), 1--59.

\pp [BH] J. Baumgartner and A. Hajnal, {\it A proof (involving Martin's
Axiom)
of a partition relation}, {\bf Fund. Math. 78} (1973),193--203.

\pp [DM] B. Dushnik and E. W. Miller, {\it Partially ordered sets}, {\bf
Amer.
J. of Math. 63} (1941), 605.

\pp [ER] P. Erd\H os and R. Rado, {\it A partition
calculus in set theory}, {\bf Bull. Amer. Math. Soc. 62} (1956), 427--489.

\pp [EHMR]  P. Erd\H os, A. Hajnal, A. M\'at\'e and R. Rado, {\bf
Combinatorial
Set Theory: Partition Relations for Cardinals}, North-Holland, 1984.

\pp [Ha] A. Hajnal, {\it Some results and problems in set theory}, {\bf Acta
Math. Acad. Sci. Hung. 11} (1960), 277--298.

\pp [Pr] K. Prikry, {\it On a problem of Erd\H os, Hajnal and Rado}, {\bf
Discrete Math. 2} (1972), 51--59.

\pp [Sh1] S. Shelah, {\it Notes on combinatorial set theory}, {\bf Israel J.
Math. 14} (1973), 262--277.

\pp [Sh2] S. Shelah, {\it On CH + $2^{\A_1}\to(\a)^2_2$ for $\a<\w_2$}, to
appear.

\bigskip
\parskip=0pt
\lineskip=0pt

Dartmouth College

Hanover, NH 03755 USA
\medskip
Mathematical Institute

Hungarian Academy of Sciences

Budapest, Hungary
\medskip
Institute Of Mathematics

Serbian Academy of Sciences

Belgrade, Yugoslavia

\bye